\documentclass[11pt]{amsart}
\usepackage[colorlinks=true,urlcolor=blue,
citecolor=red,linkcolor=blue,linktocpage,pdfpagelabels,
bookmarksnumbered,bookmarksopen]{hyperref}
\usepackage{latexsym,enumerate}
\usepackage{amsmath,amsthm,amsopn,amstext,amscd,amsfonts,amssymb}
\usepackage[latin1]{inputenc}
\usepackage[active]{srcltx}
\usepackage{verbatim}
  \newcommand{\dyle}{\displaystyle}

\newcommand{\p}{\partial}

\newcommand{\io}{\int_\Om}

\renewcommand{\a }{\alpha }

\renewcommand{\d }{\delta }
\newcommand{\e }{\varepsilon }
\newcommand{\g }{\gamma}

\newcommand{\n }{\nabla }

\newcommand{\s }{\sigma }

\renewcommand{\O }{\Omega }

\newtheorem{Theorem}{Theorem}[section]

\newtheorem{Lemma}{Lemma}[section]
\newtheorem{Proposition}{Proposition}[section]
\newtheorem{definition}{Definition}[section]
\newtheorem{remark}{Remark}[section]

\newcommand{\lap}{\Delta}

\newcommand{\car}[1]{\raise2pt\hbox{$\chi$}_{#1}}

\newcommand{\D}{\Delta }

\newcommand{\De}{\Delta}

\newcommand{\Om}{\Omega}

\long\def\salta#1{\relax}

\newcommand{\norma}[2]{\|#1\|_{\lower 4pt \hbox{$\scriptstyle #2$}}}

\newcommand{\ba}{\begin{array}}
\newcommand{\bea}{\begin{eqnarray*}}
\newcommand{\bean}{\begin{eqnarray}}
\newcommand{\be}{\begin{equation}}

\newcommand{\ea}{\end{array}}
\newcommand{\eea}{\end{eqnarray*}}
\newcommand{\eean}{\end{eqnarray}}
\newcommand{\ee}{\end{equation}}

\begin{document}
\title[Existence,   Nonexistence, and Regularity]{Existence,   Nonexistence, and Regularity Results for  Supercritical Quasilinear Elliptic Problems with Hardy Potential}
\author[B. Abdellaoui, S. Merch\'an, L. Montoro]{Boumediene Abdellaoui$^1$, Susana Merch\'an$^2$, Luigi Montoro$^3$}

\address{$^1$Boumediene Abdellaoui, Laboratoire d'Analyse Nonlin\'eaire et Math\'ematiques Appliqu\'ees, D\'epartement de Math\'ematiques,
    Universit\'e Abou Bakr Belkaid, Tlemcen 13000, Algeria.}
\email{boumediene.abdellaoui@inv.uam.es}
\address{$^2$Susana Merch\'an, Departamento de Matem\'{a}tica e Inform\'{a}tica aplicadas a las ingenier\'{i}as civil y naval, E.T.S.I. Caminos, Canales y Puertos, Universidad Polit\'{e}cnica de Madrid, Madrid, Spain.}
\email{susana.merchan@upm.es}
\address{$^3$Luigi Montoro, Dipartimento di Matematica e Informatica, Universit\`a della Calabria, Ponte Pietro Bucci 31 B, 87036 Arcavacata di Rende, Cosenza, Italy.}
\email{luigi.montoro@unical.it}



\date{}
\maketitle

\begin{abstract}
In this paper, we investigate the existence of solutions to the following supercritical elliptic problem involving the $p$-Laplace operator  and the Hardy-Leray potential
\begin{equation*}
\left\{\begin{array}{rcll}
    -\Delta_p u +u^\beta|\nabla u|^p &= & {\vartheta}  \dfrac{u^q }{|x|^p}+ f\;\; &\hbox{in } \Omega\\
    u &\geq & 0\;\; &\hbox{in } \Omega \\
    u &= & 0\;\;&\hbox{on } \partial\Omega,
\end{array}\right.
\end{equation*} where  $0\in\Omega$. We analyze the effect of the perturbation term $u^\beta|\nabla u|^p$, with $\beta \geq 0$, and the role of the  exponent $q$  on the existence of solutions. We also analyze nonexistence of solutions  when $f$ is a singular measure.
Finally, we  focus on the regularity and asymptotic  properties of these solutions.
\end{abstract}
\section{Introduction}
\noindent In this paper, we  focus on the existence and regularity  of  positive weak  solutions (see Definition \ref{def:weak} below) to the supercritical problem
\begin{equation}\tag{$\mathcal P$}\label{te}
\left\{\begin{array}{rcll}
    -\Delta_p u +u^\beta|\nabla u|^p &= & {\vartheta}  \dfrac{u^q }{|x|^p}+ f\;\; &\hbox{in } \Omega\\
    u &\geq & 0\;\; &\hbox{in } \Omega \\
    u &= & 0\;\;&\hbox{on } \partial\Omega,
\end{array}\right.
\end{equation}
where $\Omega$ is a bounded and smooth domain, $1<p<N$, $\beta\geq 0$, $\vartheta>0$,  $f\in L^1(\Omega)$, $f\gneqq0$ in~$\Omega$ and
$$0\in\Omega.$$
The role of the perturbation term in this type of problem is fundamental, as the pure supercritical problem, i.e. when $q>p-1$ and the perturbation term $u^\beta|\nabla u|^p $ is absent,  does  not admit solutions, not even in the weakest sense such as entropy solutions \cite{BBGGPV}. Indeed,  if we consider the approximating problem associated with $\eqref{te}$,
\begin{equation*}\label{Ireneopnn}
\left\{\begin{array}{rcll}-\Delta_p u_n &= & \dfrac{u^q_n }{(1+\frac 1n )|x|^p}+f_n\;\;&\hbox{in } \Omega\\
 u_n &\geq & 0\;\;&\hbox{in } \Omega \\
 u_n &= & 0\;\;&\hbox{on }\partial\Omega,
\end{array}\right.
\end{equation*}
where $q>p-1$,  $0\leq f_n\in L^\infty(\O)$ and $f_n\rightarrow f \in L^1(\Omega)$, then as in
\cite{Abdellaoui-Peral-20031}, using a variant of the weak Harnack
inequality, it can be proved that  $u_n(x_0)\to \infty$ as $n\to \infty$
for every $x_0\in \O$.

Despite this last fact, it is important to point out  that the pure supercritical  problem (i.e.  problem \eqref{te} with  $q>p-1$ and without the perturbation term $u^\beta|\nabla u|^p $)
 exhibits a very different  behaviour when $0\in\partial\Omega$.\\
 In fact, under suitable geometric  assumptions on the   domain $\Omega$,  the existence of energy  solutions can be established, see   \cite{DPe, MMP}.

\

On the other hand, the situation in the subcritical case is more standard, as solutions can be obtained  either by classical methods
based on the  minimization of the associated energy functional,  or  by approximation techniques; see also $(ii)$ of Remark \ref{rem:min}.

\

{Related to the quasilinear case, problems
involving natural growth in the gradient arise in a variety of
physical models, including non-Newtonian fluid dynamics, heat
conduction, and reaction-diffusion systems with convection
effects. Such terms typically of order $|\nabla u|^q$, with $q$
matching the principal elliptic operator introduce significant
analytical challenges. In particular, the lack of coercivity and
the failure of standard energy estimates for these nonlinearities
make the existence and regularity theory highly nontrivial.}

{A foundational contribution by Ferone and Murat~\cite{FM} established existence results under natural growth
assumptions when the source term remains suitably small. Boccardo~\cite{B00} examined model equations of the form
$$-\text{div}((1+|u|^r)\n u)+|u|^{m-2}|\n u|^2=f,$$
and derived existence results in appropriate function spaces. The quasilinear counterpart of such problems,
with the inclusion of singular weights, was investigated in
\cite{ABPP0}, where global existence results were established for
a broad class of quasilinear elliptic equations. By contrast, in
the sub-natural growth case (i.e., when $0<q<p-1$), the problems
are better posed and can often be treated using more classical
tools. Together, these works illustrate a nuanced landscape:
natural-gradient growth requires delicate approximation and
compactness methods, whereas sub-natural growth admits a more
standard variational treatment.}

{When a Hardy potential is present typically on the
form $\lambda\frac{u^a}{|x|^\theta}$ an additional layer of
difficulty arises. By the Hardy inequality, such singular
terms may prevent existence or trigger blow-up phenomena for large
$\lambda$ or critical values of $a$. Thus, the Hardy potential
often acts as a strong obstruction to solvability and regularity.}

However, a mechanism that can
neutralize this singularity term from gradient absorption terms.
Nonlinearities of the type $u^\beta|\nabla u|^p$ can act as an
absorption effect, counterbalancing the Hardy potential and
restoring coercivity and compactness properties. The interaction
between the singular potential and the absorption term induced by the
gradient is, therefore, decisive in determining the existence of weak
solutions exist. For the particular case $\beta=0$, relevant
results are available in \cite{APP}, \cite{CMMS},  \cite{MMP} and
\cite{simmetry}.

{In this framework, our contribution focuses on the
combined effect of a Hardy-singular source and a natural-growth
convection term of the form $u^\beta|\nabla u|^p$. We establish existence, regularity, and stability
results within parameter ranges for $(\beta,q,\vartheta)$,
under minimal assumptions on the data $f$. We also highlight
parameter regimes where the Hardy constant plays a critical role
in solvability. Moreover, we analyze  nonexistence phenomena  when the datum
$f$ is a singular measure.}

\

In order to establish the framework of our paper, we begin by
introducing the definition of weak solution.
\begin{definition}\label{def:weak}
    We say that a nonnegative measurable function $u$ is a weak solution to~\eqref{te},
    if $u\in W^{1,p}_0(\Omega)$,
    \[u^\beta|\nabla u|^p\in L^1(\Omega),\qquad \dfrac{u^q }{|x|^p}\in L^1(\Omega)\] and
    $$
    \int_{\Omega}|\nabla u|^{p-2}(\nabla u,\nabla \phi)+\int_{\Omega}u^\beta|\nabla u|^p\phi=\vartheta \int_{\Omega}\dfrac{u^q }{|x|^p}\phi+ \int_{\Omega}f\phi\quad \forall \phi \in W^{1,p}_0(\Omega)\cap L^{\infty}(\Omega).$$
\end{definition}
The following theorem provides the first main result of our paper.
\begin{Theorem}\label{t1}
 Consider the problem \eqref{te} with $1<p<N$, $\beta \geq 0$,  and assume
 \[ f\in L^1(\Omega) \quad \text{and}\quad f \gneqq0 \,\,\text{in}\,\, \Omega.\]
\begin{itemize}
\item [($\ast$)] If $p-1 < q<p+\beta$, then, for all $\vartheta>0$, there exists a
 positive weak solution $u\in W_0^{1,p}(\Omega)$ to \eqref{te}.

\

\item[($\ast\ast$)] If $q=p+\beta$, then, there exists $\bar\vartheta=\bar\vartheta(p,\beta)$ such that,  for all $0<\vartheta< \bar\vartheta$,
there exists a   positive weak solution $u\in W_0^{1,p}(\Omega)$ to \eqref{te}.
\end{itemize}
\end{Theorem}
\begin{remark}\label{rem:min}We point out the following:
\begin{itemize}
\item [$(i)$]The solution $u$ obtained in Theorem \ref{t1} is the minimal solution of problem~\eqref{te},  in the sense that  any other solution or supersolution $\hat{u}$ of \eqref{te} satisfies
$u\leq\hat{u}$ in~$\Omega$. This will become clear in the proof of Theorem \ref{t1},  by the very definition of the approximating sequence for $u$, see \eqref{eq:ta}.
\item [$(ii)$] Theorem \ref{t1} holds for every
\(0 \leq q \leq p + \beta\). In this paper, however, we restrict to the range
\(p - 1 < q \leq p + \beta\), namely the supercritical case, since it is both more
interesting and challenging. In fact, in the range \(0 \leq q \leq p - 1\),
the existence of solutions can be obtained by standard techniques.
\end{itemize}
\end{remark}
In order to prove Theorem \ref{t1}, we first establish the existence of solutions to a related problem similar to \eqref{te}, but without the Hardy singular term, see Theorem \ref{thone} below.

{The situation changes dramatically, however, when the datum
$f$ is a singular measure with respect to the $p$-capacity
$\operatorname{cap}_p$, (see  Definition \ref{pcap} and Definition \ref{singular} below). In this setting, even in the absence of the
Hardy singular term, one can establish a nonexistence result: there are  no {\em reasonable} solutions in a  strong sense. Namely, if we approximate $f$ by smooth
functions in the narrow topology of measures, then the corresponding approximate solutions
converge to 0. We prove a result of the following type:}
\begin{Theorem}\label{sing} Let $\mu$ be a nonnegative singular Radon measure, that is,   $\mu$ is concentrated on a Borel
set {$E\subset \subset \O$},  such that  $\operatorname{cap}_{p}(E) = 0$. Let $\{f_n\}_n\subset L^1(\O)$ be a non negative sequence  such that $f_n\rightharpoonup \mu$ in the narrow topology of measures.
Let $u_n$ denotes   the unique solution to the problem
\begin{equation}\label{apr01pp}
\left\{\begin{array}{rcll}
-\De_p u_n+u^\beta_n|\nabla u_n|^p&=&f_n\;\; &\hbox{in } \Omega\\
u_n &> & 0\;\; &\hbox{in } \Omega \\
u_n &= & 0\;\;&\hbox{on } \partial\Omega.\end{array}\right.
\end{equation}
 Then $u_n\rightharpoonup 0$ weakly
in $W_0^{1,p}(\Omega)$.
\end{Theorem}
The situation described in Theorem \ref{sing}  also extends to  the general problem \eqref{te}. We analyze this case in   Theorem \ref{singthh} here below. This result can be regarded  as a nonexistence results for \eqref{te} when $f$ is a singular measure, in the sense that suitable approximating solutions does not converge to a {\em rasonable} limit, namely the minimal solution of \eqref{te} with $f=0$. In this case, the problem  also admits  the trivial solution.
\begin{Theorem}\label{singthh} Let $\mu$ be a nonnegative singular Radon measure, that is,   $\mu$ is concentrated on a Borel
set {$E\subset \subset \O$},  such that  $\operatorname{cap}_{p}(E) = 0$. Let $\{f_n\}_n\subset L^1(\O)$ be a non negative sequence  such that $f_n\rightharpoonup \mu$ in the narrow topology of measures.
Let $u_n$ denote the minimal solution to the
problem
\begin{equation}\label{sin0111}
\left\{\begin{array}{rcll}
-\De_p u_n+u^\beta_n|\nabla u_n|^p&=& \vartheta\dfrac{u_n^q}{|x|^p}+f_n\;\; &\hbox{in } \Omega\\
u_n &> & 0\;\; &\hbox{in } \Omega \\
u_n &= & 0\;\;&\hbox{on } \partial\Omega,\end{array}\right.
\end{equation}
as obtained in Theorem \ref{t1}. Then $u_n\rightharpoonup \hat{u}$
weakly in $W_0^{1,p}(\Omega)$, where $\hat{u}$ denotes a solution (possibly  trivial) of
problem~\eqref{te} with $f=0$.
\end{Theorem}
\

In the second part of this work, we  focus on the regularity and asymptotic behavior of solutions. The main results are the following:
\begin{Theorem}\label{l1}
Let $u$ be a non-negative weak solution  to problem \eqref{te}.  Then,
    $$ \lim_{|x|\rightarrow 0}u(x)= +\infty.$$
\end{Theorem}
Theorem \ref{l1}  shows that, due to the presence of the Hardy potential, every solution (in fact, every supersolution, as one can see in the proof) is unbounded, even in the case of regular data
$f$. This stands in contrast with the classical regularity theory for problems without the Hardy potential.

The next theorem is a Calder\'on-Zygmund type result for the minimal solution to problem \eqref{te}.
\begin{Theorem}\label{thm2bis}
 Assume that the conditions of Theorem \ref{t1} hold with  \[ f\in L^m(\Omega),\,\, m>1 \quad \text{and}\quad f \gneqq0 \,\,\text{in}\,\, \Omega.\]
Let $u$  be the minimal solution of problem \eqref{te}.

\begin{itemize}
\item In  case $(*)$, the following assertions hold:
\begin{itemize}
\item[$(i)$] If $1<m<\frac{N}{p}$, then $u\in L^\g(\Omega)$ for all $\gamma\le \frac{mp^*(\beta+p)}{pm-p^*(m-1)}$. Moreover,
    $$\io u^{\frac{p(m\beta+p^*(m-1))}{pm-p^*(m-1)}}|\n u|^p dx<\infty.$$
   \item[$(ii)$] If $m\ge \frac{N}{p}$, then $u\in L^\gamma(\Omega)$ for all $\gamma<\infty$. Moreover,
    $$\io u^{\gamma}|\n u|^p dx<\infty,\quad\mbox{  for all  }\gamma<\infty.$$
    \end{itemize}
 \item In  case $(**)$, there exists $\bar{\bar\vartheta}<{\bar\vartheta}$  such that the statements  $(i)$ and $(ii)$ remain valid  for all $0<\vartheta< \bar{\bar\vartheta}$.
 \end{itemize}
\end{Theorem}
The article is organized as follows:
\begin{itemize}
    \item In Section \ref{sec01} we prove our main existence result, namely Theorem \ref{t1}. This section begins with the study of existence and nonexistence of solutions of a variant of problem \eqref{te} (see problem \eqref{one} below). The preliminary existence result obtained here is used to establish the main existence theorem for the general problem \eqref{te}. We also provide  the proof of the nonexistence results stated in Theorem \ref{sing} and Theorem \ref{singthh}.
\item In Section \ref{sec3}, we study certain regularity properties of the solutions and their  behaviour near the singular point zero. In this context,  we prove Theorem~\ref{l1} and Theorem \ref{thm2bis}.  Moreover, we provide precise asymptotic descriptions of the solutions for specific choices of the data $f$.
\end{itemize}

\

\noindent{\sc Notation.} Generic fixed numerical constants will be denoted by $C$ (sometimes with subscript) and may vary from line to line or within a single formula.  Moreover $f^+$ and $f^-$ denote the positive and negative parts of a function $f$, respectively; that is
$f^+=\max\{f,0\}$ and $f^-=\min\{f,0\}.$ We also denote by $|A|$ the Lebesgue measure of the set $A$.
\section{Existence of Solutions to \eqref{te}}\label{sec01}
{In this section, we begin by recalling and presenting some technical results, together with a few compactness results.

We first introduce   the following truncation function, which will be used systematically throughout all the paper:
\begin{equation}\label{deftroncata}
T_k(s)=\max\{\min\{k,s\},-k\}, \quad k>0.
\end{equation}
The first compactness result we use   is established in Theorem 3.4 of \cite{DMOP}.
\begin{Theorem}\label{compaggth}
Let $\{f_n\}_n\subset L^\infty(\O)$ be such that
$||f_n||_{L^1(\O)}\le C$ for all $n$. Define $u_n$ to be the unique weak
solution to the problem
\begin{equation*}\left\{
\begin{array}{rcll}
    -\Delta_p u_n &=& f_n & \mbox{ in } \Omega\\
    u_n &= & 0 & \mbox{ on } \partial\Omega.
\end{array}
\right.
\end{equation*}
Then \begin{enumerate} \item the sequence $\{|u_n|^{p-2}u_n\}_n$
is bounded in the space $L^{\sigma}(\O)$ for all
$\sigma<\frac{N}{N-p}$; \item the sequence $\{|\n u_n|^{p-2}\n
u_n\}_n$ is bounded in $L^\theta(\O)$ for all
$\theta<\dfrac{N}{N-1}$; \item for all $k>0$, the sequence
$\{T_k(u_n)\}_n$ is bounded in $W^{1,p}_0(\O)$.
\end{enumerate}
Moreover, if, up to a subsequence $f_n\rightharpoonup f\in
L^1(\O)$, there
exists a measurable function $u$ such that $T_m(u)\in W^{1,p}_0(\O)$ for every $m>0$ and,  up to subsequence, we have $T_m(u_n)\to T_m(u)$ strongly
in $W^{1,p}_0(\O)$ for every $m>0$.
\end{Theorem}
}
{The second  compactness result that will be used
systematically in this paper concerns to the strong convergence
of the truncating sequences for increasing sequences. The case $p=2$ was proved in \cite{AP} and \cite{Murat1}. We  include here the proof for the nonlinear case $p\neq 2$.
\begin{Theorem}\label{APM}
Let $\{w_n\}_n\subset W^{1,p}_0(\O)$ be a sequence of nonnegative
increasing functions such that $-\D_p w_n=g_n\ge 0$ in $\O$.
Suppose that the sequence $\{T_k(w_n)\}_n$ is bounded in
$W^{1,p}_0(\O)$ for all $k>0$. Then, there exists a measurable
function $w$ such that $T_k(w)\in W^{1,p}_0(\O)$ and
$T_k(w_n)\rightarrow T_k(w)$ strongly in $W^{1,p}_0(\O)$, for
every $k>0$.
\end{Theorem}
\begin{proof}
Since $\{w_n\}_n$ is an increasing sequence, then, using the fact
that $\{T_k(w_n)\}_n$ is bounded in $W^{1,p}_0(\O)$ for all $k>0$,
we get the existence of a measurable function $w\ge 0$ such that
$T_k(w)\in W^{1,p}_0(\O)$ for all $k>0$, $T_k(w_n)\le T_k(w)$ and,
up to a subsequence, $T_k(w_n)\rightharpoonup T_k(w)$ weakly in
$W^{1,p}_0(\O)$. Since the weak limit is unique,
the whole sequence $\{T_k(w_n)\}_n$ converges weakly to $T_k(w)$.
Now, we prove the strong convergence. Notice that
$$ \|\n T_k(w)\|_{L^p(\O)}\le \liminf_{n\to \infty}\|\n
T_k(w_n)\|_{L^p(\O)}, \mbox{ for all }k>0.
$$
Now, since $T_k(w_n)\le T_k(w)$, we deduce that
$$
 \io -\D_p w_n
(T_k(w_n)) dx \le \io -\D_p w_n T_k(w) dx.
$$
Thus
$$
\io |\n T_k(w_n)|^p dx\le \io |\n T_k(w_n)|^{p-1}|\n T_k(w)| dx
\le \dyle \frac{1}{p'}\io |\n T_k(w_n)|^pdx + \frac{1}{p}\io |\n
T_k(w)|^pdx. $$ Then
$$
\io |\n T_k(w_n)|^pdx \le\io |\n T_k(w)|^pdx.
$$
Hence, we conclude that
$$
\limsup\limits_{n\to \infty}\io|\n T_k(w_n)|^p\le \io |\n
T_k(w)|^pdx.
$$
Therefore, we conclude that
$$
\limsup\limits_{n\to \infty}\io|\n T_k(w_n)|^p dx
=\liminf\limits_{n\to \infty}\io|\n T_k(w_n)|^p dx =\io |\n
T_k(w)|^p dx.
$$
{Since $W^{1,p}_0(\O)$ is a uniformly convex space},
it follows that $T_k(w_n)\to T_k(w)$ strongly in $W^{1,p}_0(\O)$.
\end{proof}}
We recall a useful comparison principles that will be used in what follows.
\begin{Proposition}[\cite{Abdellaoui-Peral-20031}]\label{rome2}
{Assume that $w_1,w_2\in W^{1,p}(\Omega)$ are non negative functions such that $w_1\ge w_2$ on $\p\Omega$,
\[\frac{w_1}{w_2},\,\,\frac{w_2}{w_1}\in L^\infty(\Omega)\]
and
$$
\left\{
\begin{array}{lll}
-\D_p w_1 &\ge & h_1(x)K(w_1) \mbox{  in }\O\\
-\D_p w_2 &\le & h_2(x)K(w_2) \mbox{  in }\O,
\end{array}\right.
$$
where $h_1,h_2\in L^{\infty}(\Omega)$ are non negative  functions satisfying $h_1\geq h_2$ in $\Omega$ and $h_iK(w_i)\in W^{-1,p'}$, $i=1,2$. Suppose further that $K: \mathbb{R}^+\to \mathbb{R}^+$ is  continuous and that the map
\[\sigma\rightarrow \frac{K(\sigma)}{\sigma^{p-1}}\] is strictly decreasing.
Then $w_1\ge w_2$ in $\Omega$.}
\end{Proposition}
We prove a useful  comparison principle.
\begin{Proposition}\label{compagg}
Suppose that $0\leq \theta<p-1$ and $\beta\geq0$. Let $u_1,u_2\in W^{1,p}(\Omega)$ are non negative functions such that $u_1\ge u_2$ on $\p\Omega$, and assume that
\[u_1^\beta|\n u_1|^p, \,\,u_2^\beta|\n u_2|^p\in L^1(\Omega)\]
and
    $$
    \left\{
    \begin{array}{lll}
        -\D_p u_1 +u_1^\beta|\n u_1|^p &\ge & h_1(x)u_1^\theta\mbox{  in }\O\\
        -\D_p u_2 +u_2^\beta|\n u_2|^p &\le & h_2(x)u_2^\theta\mbox{  in }\O,
    \end{array}\right.
    $$
where $h_1,h_2\in L^{\infty}(\Omega)$ are non negative  functions satisfying $h_1\geq h_2$ in $\Omega$.    Then, $u_1\ge u_2$ in $\Omega$.
\end{Proposition}
\begin{proof}
Assuming the above conditions, for $i=1,2$, we define $
v_i=H(u_i)$ where $H:\mathbb R^+_0\rightarrow \mathbb R$, is defined by
$$
H(\sigma)=\int_0^\sigma
e^{-\frac{t^{\beta+1}}{(p-1)(\beta+1)}}\, dt.
$$
It is clear that $H$ is strictly increasing concave function.
Setting $D(s)=H^{-1}(s)$, then $v_i\in L^\infty(\O)\cap W^{1,p}(\O)$ and $v_i$ satisfies
    $$
\left\{
\begin{array}{lll}
    -\D_p v_1&\ge & h_1(x)\frac{D(v_1)^\theta}{(D'(v_1))^{p-1}}\mbox{  in }\O\\
    -\D_p v_2&\leq & h_2(x)\frac{D(v_2)^\theta}{(D'(v_2))^{p-1}}\mbox{  in }\O,
\end{array}\right.
$$
with $v_1\ge v_2$ on $\p\Omega$. For $s>0$, let
\[R(s)=\frac{D(s)^\theta}{s^{p-1}(D'(s)^{p-1}},\] then we claim that $R$ is decreasing. Notice that $R(s)=\tilde R(\sigma (s))$, where
$$
\tilde{R}(\sigma)=\dfrac{\sigma^\theta(H'(\sigma))^{p-1}}{(H(\sigma))^{p-1}}.
$$
Thus $R'(s)=\tilde R'(\sigma (s))\sigma '(s)= \tilde R'(\sigma (s))\ D'(s).$
Taking into consideration that $D'(s)>0$, then we only have to claim that $$\tilde{R}'(\sigma (s)) < 0.$$  Notice that
$$\tilde{R}(\sigma)=\sigma^{\theta-(p-1)}\bigg(\dfrac{\sigma H'(\sigma)}{H(\sigma)}\bigg)^{p-1}.$$
Since $\theta<p-1$, to show the claim we only have to prove that the derivative of $\dfrac{\sigma H'(\sigma)}{H(\sigma)}$ is  negative. Hence, we have to consider the sign of
$$
I(\sigma)=H(\sigma)\bigg(H'\sigma+\sigma H''(\sigma)\bigg)-\sigma(H'(\sigma))^2.
$$
Since $H''(\sigma)=-\frac{\sigma^\beta}{p-1} H'(\sigma)$, the sign of $I$
is the sign of
$$J(\sigma)=H(\sigma)\left(1-\frac{\sigma^{\beta+1}}{p-1}\right)-\sigma H'(\sigma).$$
After some computations, we obtain that   $J'(s)=-\frac{(\beta+1)}{p-1}H(\sigma)\sigma^\beta\le 0$.
Since $J(0)=0$, then $J(\sigma)< 0$ for $\sigma > 0$ and then the claim
follows.

Thus, $R$ is (strictly) decreasing. Therefore, using the comparison
principle in Proposition~\ref{rome2}, we deduce that $v_1\ge v_2$
and, then, $u_1\ge u_2$. Hence, we conclude.
\end{proof}
Finally, in the next lemma, we recall some useful algebraic inequalities related to the
$p$-Laplacian operator (see \cite{lucio}), which will be used throughout the paper.
\begin{Lemma}\label{HHH}
There exists a constant $C(p) >0$ such that for all $\xi_1,\xi_2\in \mathbb{R}^N$, with $|\xi_1|+|\xi_2|>0$
such that
\begin{itemize}
\item $
\langle |\xi_1|^{p-2}\xi_1-|\xi_2|^{p-2}\xi_2,
\xi_1-\xi_2\rangle\ge
C(p){(|\xi_1|+|\xi_2|)^{p-2}}{|\xi_1-\xi_2|^2}$, if $1<p<2$;

\

\

\item $\langle |\xi_1|^{p-2}\xi_1-|\xi_2|^{p-2}\xi_2,
\xi_1-\xi_2\rangle\ge C(p)|\xi_1-\xi_2|^p$,  if $p\geq 2$.
\end{itemize}
\end{Lemma}
\subsection{Problem without the Hardy potential: Existence and Nonexistence Results}
In this section  we consider problem \eqref{one} (stated below), which corresponds to problem  \eqref{te} without the Hardy potential. We first prove the existence and uniqueness  of a solution to  \eqref{one}. This result, besides being interesting in its own right,  provides a compactness criterion that will be essential in the next section to show the existence of solutions to problem \eqref{te}, which will in turn lead to the proof of our main existence result, Theorem \ref{t1}.

Finally, at the end of this section, we prove the nonexistence result stated in Theorem~\ref{sing}, in the case where the $L^1$ datum $f$ is replaced by a singular measure $\mu$.

\

Let us consider the following problem
 \begin{equation}\label{one}
\left\{\begin{array}{rcll}
    -\Delta_p u +u^\beta|\nabla u|^p &=& f\;\; &\hbox{in } \Omega\\
    u &\geq & 0\;\; &\hbox{in } \Omega \\
    u &= & 0\;\;&\hbox{on } \partial\Omega,
\end{array}\right.
\end{equation}
where {{$1<p<N$}, $\beta\geq 0$ and $f\in L^1(\Omega)$ with
$f\gneqq0$.
The existence result is stated below. We would like to point out that Theorem \ref{thone}, as can be verified in its proof, actually holds for all $p > 1$. Nevertheless, since this result will be used in the proof of Theorem \ref{t1}, we prefer to state it under the same assumptions as the main problem \eqref{te} involving the Hardy potential, where the condition $1 < p < N$ is necessary. We adopt this restriction purely for the  reader's convenience.
\begin{Theorem}\label{thone}
Assume $1<p<N$, $\beta\geq 0$ and $f\in L^1(\Omega)$ with
$f\gneqq0$. Then, problem~\eqref{one}
has a unique positive solution   $u\in W^{1,p}_0(\O)$ such that
\[u^\beta|\n u|^p \in L^1(\O).\] Moreover, if $0\lneqq f_1\le f_2$,
then $u_1\le u_2$.

{Furthermore, if $\{f_n\}_n \subset L^\infty(\O)$} is
an increasing sequence such that $f_n\to f$  in $L^1(\O)$,
denoting $u_n$ (resp. $u$) the solution of problem \eqref{one}
with datum $f_n$ (resp. $f$), then the sequence $\{u_n\}_n$ is  monotone increasing and
\[u_n\to u  \,\,\text{in }
W^{1,p}_0(\O) \qquad \text{and} \qquad u_n^\beta|\n u_n|^p\to u^\beta|\n u|^p\,\,
\text{in } L^1(\O).\]
\end{Theorem}
\begin{proof}
To prove the
existence we follow by approximation and we give the proof in two steps.

\

\begin{enumerate}
\item {\emph{ Existence and uniqueness for bounded datum.}} First
of all let us assume that $f\in L^\infty(\O)$ is a nonnegative
function and $u_k$ the solution of the problem
\begin{equation}\label{apr0}
\begin{cases}
-\De_p u_k+\dfrac{(u_k^+)^{\beta}|\nabla
u_k|^p}{1+\frac{1}{k}(u_k^+)^{\beta}|\nabla u_k|^p}=f& \text{in }
\Omega
\\
u_k=0 & \text{on } \partial \Omega.
\end{cases}
\end{equation}
Notice that the existence of $u_k$, for $\beta>0$  follows using
the classical Leray-Lions theorem \cite{LeLi}. Since  $f\gneq 0$,
and using $u^-_k$ as a test function in \eqref{apr0}, it follows
that $u_k\ge 0$.

{{For the case $\beta=0$, we can not use  the
Leray-Lions theorem since the non-linear term is not continuous in
$\n u$}}. To show the existence in this case we use properly  Schauder fixed
theorem.

Define the operator $\mathcal{T}: W^{1,p}_0(\O)\to W^{1,p}_0(\O)$
by setting $u=\mathcal{T}(v)$ and $u$ the (unique) solution to
\begin{equation*}
\left\{
\begin{array}{rcll}
-\De_p u &= & f-\dfrac{|\nabla v^+|^p}{1+\frac{1}{k}|\nabla v^+|^p}&
\text{in } \Omega
\\
u &= & 0 & \text{on } \partial \Omega.
\end{array}
\right.
\end{equation*}
The existence of $u$ follows using a simple minimizing argument.
Notice that
$$
||u||_{W^{1,p}_0(\O)}\le C(\O)\bigg(\io (k+|f|)^{p'}
dx\bigg)^{\frac{1}{p}}.
$$
Thus, setting $R=C(\O)\bigg(\io (k+|f|)^{p'}
dx\bigg)^{\frac{1}{p}}$, (being $p'=p/(p-1)$) it follows that
$\mathcal{T}(\overline{B_R(0)})\subset \overline{B_R(0)}$.

\begin{itemize}
\item {We prove that $\mathcal{T}$ is continuous.}
Let $\{v_n\}_n\subset W^{1,p}_0(\O)$ such that $v_n\to v$ strongly
in $W^{1,p}_0(\O)$. Notice that, by Theorem 1, in \cite{MM}, we
have that $v_n^+\to v^+$ strongly in $W^{1,p}_0(\O)$. Thus, by
the dominated convergence theorem we deduce that
$$
\dfrac{|\nabla v_n^+|^p}{1+\frac{1}{k}|\nabla v_n^+|^p}\to
\dfrac{|\nabla v^+|^p}{1+\frac{1}{k}|\nabla v^+|^p}
$$
strongly in {$L^{p'}(\O)$}. Let us define $u_n=\mathcal{T}(v_n)$,
$u=\mathcal{T}(v)$ and set
\begin{equation}\label{eq:sire}
g_n=-\dfrac{|\nabla v_n^+|^p}{1+\frac{1}{k}|\nabla
v_n^+|^p}+f,
\end{equation}
then $g_n\to g$ strongly in ${L^{p'}(\O)}$, where
$g=-\dfrac{|\nabla v^+|^p}{1+\frac{1}{k}|\nabla v^+|^p}+f$.
Hence, we conclude that $u_n\to u$ strongly in $W^{1,p}_0(\O)$.
Therefore, $\mathcal{T}$ is continuous.

\item {We prove  that $\mathcal{T}$ is compact.} Let
$\{v_n\}_n\subset W^{1,p}_0(\O)$ be a bounded sequence, then, there
exists $C>0$ such that $||v_n||_{W^{1,p}_0(\O)}\le C$.
{Taking into consideration that $||g_n||_{{L^{p}(\O)}}\le C(k,f)$, see \eqref{eq:sire}, we have that up to a subsequence $g_n \rightharpoonup g$ in $L^p(\Omega)$. Then
$$\int_{\Om} g_n \varphi \,dx \ \xrightarrow{n \rightarrow +\infty } \int_{\Om}g\varphi \,dx \quad \forall \varphi \in L^{p'}(\Omega),$$
and therefore
$$\int_{\Om} g_n \varphi \,dx \ \xrightarrow{n \rightarrow +\infty } \int_{\Om}g\varphi \,dx \quad \forall \varphi \in L^{\infty}(\Omega),$$
namely
$g_n \rightharpoonup g$ in $L^1(\Omega)$. }
Using Theorem \ref{compaggth}, we get the existence of $u\in
W^{1,p}_0(\O)$ such that, for all $i>0$, we have $T_i(u_n)\to
T_i(u)$ strongly in~$W^{1,p}_0$.
Since $\|u_n\|_{L^{\infty}(\Omega)}\le
C(k,||f||_{L^\infty(\O)})$ (using for instance a Stampacchia-type argument) and choosing
$i>C(k,||f||_{L^\infty(\O)})$, it follows that $u_n\to u$ strongly
in $W^{1,p}_0(\O)$. Hence, $\mathcal{T}$ is compact.
\end{itemize}
Therefore, using the Schauder fixed point Theorem we deduce the
existence of $u_k\in \overline{B_R(0)}$ such that
$\mathcal{T}(u_k)=u_k$. Thus, $u_k$ solves
\begin{equation}\label{apr0fix2}
\begin{cases}
-\De_p u_k+\dfrac{|\nabla u^+_k|^p}{1+\frac{1}{k}|\nabla u^+_k|^p}=f&
\text{in } \Omega
\\
u=0 & \text{on } \partial \Omega.
\end{cases}
\end{equation}
Since $f\gneq 0$ and using $u^-_k$ as a test function in
\eqref{apr0fix2}, we deduce that $u_k\ge 0$.

Hence, we conclude that for all $\beta\geq 0$, the problem

\begin{equation}\label{LLL}
\begin{cases}
-\De_p u_k+\dfrac{u_k^\beta|\nabla
u_k|^p}{1+\frac{1}{k}u_k^{\beta}|\nabla u_k|^p}=f& \text{in }
\Omega
\\
u_k\ge 0 & \text{in } \partial \Omega\\
u_k=0 & \text{on }
\partial \Omega,
\end{cases}
\end{equation}
has a nonnegative solution such that $u_k\in W^{1,p}_0(\O)$.

By comparison, it is clear that $u_k\le w$, where $w$ is the unique
bounded solution to the problem
$$
-\D_p w=f,\,\,\, w\in W^{1,p}_0(\O)\cap L^\infty(\O).
$$
Using $u_k$ as a test function in \eqref{LLL}, it follows that
$$
\io |\n u_k|^p dx\le \io f u_k dx\le
C(\Omega)||w||_{L^\infty(\O)}||f||_{L^\infty(\O)}.
$$
{Hence, we deduce that the sequence $\{u_k\}_k$ is bounded in
$W^{1,p}_0(\O)\cap L^\infty(\O)$. As a result, up to a subsequence, there exists $u\in W^{1,p}_0(\O)\cap L^\infty(\O)$ such
that $u_k\rightharpoonup u$ weakly in $W^{1,p}_0(\O)$ and $u_k\to
u$ strongly in $L^r(\O)$ for all $r>1$.}

We claim that $u_k\to u$ strongly in $W^{1,p}_0(\O)$. To prove the
claim we consider first $(u_k-u)^+$ as a test function in \eqref{LLL},
then
$$
\io |\n u_k|^{p-2}(\n u_k,\n (u_k-u)^+ )\, dx \le \io f(u_k-u)^+ dx.
$$
Using the fact that, by weak convergence
$$
\io |\n u|^{p-2}(\n u,\n (u_k-u)^+)\, dx=o(1) \qquad\mbox{and} \qquad\io
f(u_k-u)^+ dx=o(1),
$$
we deduce that
$$
\io (|\n u_k|^{p-2}\n u_k- |\n u|^{p-2}\n u)\n (u_k-u)^+\, dx
\le o(1).
$$
Thus, by applying Lemma \ref{HHH},  we finally obtain
\[\io |\n
(u_k-u)^+|^p \le o(1).\] Therefore,  $(u_k-u)^+\to 0$ strongly in
$W^{1,p}_0(\O)$.

Related to $(u_k-u)^-$, using $e^{-u^{\beta+1}_k}(u_k-u)^-$ as a
test function in \eqref{LLL}, we obtain
$$
\begin{array}{lll}
& \dyle \io e^{-u^{\beta+1}_k}(|\n u_k|^{p-2}\n u_k,\n (u_k-u)^-)\,
dx\\
& +\dyle \io e^{-u^{\beta+1}_k}\bigg(\dfrac{u_k^\beta|\nabla
u_k|^p}{1+\frac{1}{k}u_k^{\beta}|\nabla
u_k|^p}-(\beta+1)u_k^\beta|\nabla u_k|^p\bigg)(u_k-u)^- \, dx\\
&=\dyle \io f e^{-u^{\beta+1}_k}(u_k-u)^- dx.
\end{array}
$$
Since $(u_k-u)^-=-(u-u_k)^+$, it follows that
$$
\begin{array}{lll}
&\dyle \io e^{-u^{\beta+1}_k}|\n u_k|^{p-2}(\n u_k,\n (u-u_k)^+)
\, dx\\
& +\dyle \io e^{-u^{\beta+1}_k}\bigg((\beta+1)u_k^\beta|\nabla
u_k|^p-\dfrac{u_k^\beta|\nabla
u_k|^p}{1+\frac{1}{k}u_k^{\beta}|\nabla
u_k|^p}\bigg)(u_k-u)^- \, dx\\
&=-\dyle \io f e^{-u^{\beta+1}_k}(u_k-u)^- \, dx.
\end{array}
$$
Notice that
$$
\dyle \io e^{-u^{\beta+1}_k}\bigg((\beta+1)u_k^\beta|\nabla
u_k|^p-\dfrac{u_k^\beta|\nabla
u_k|^p}{1+\frac{1}{k}u_k^{\beta}|\nabla u_k|^p}\bigg)(u_k-u)^-
\, dx\leq 0,
$$
then
$$
\dyle \io e^{-u^{\beta+1}_k}|\n u_k|^{p-2}(\n u_k,\n (u-u_k)^+)\,
dx\geq -\dyle \io f e^{-u^{\beta+1}_k}(u_k-u)^-\, dx=o(1).
$$
Using again the fact that
$$
\dyle \io e^{-u^{\beta+1}_k}|\n u|^{p-2}(\n u,\n (u-u_k)^+)\,  dx=o(1),
$$
we reach that
$$
\dyle \io e^{-u^{\beta+1}_k}( |\n u|^{p-2}\n u-|\n
u_k|^{p-2}\n u_k)\n (u-u_k)^+ \, dx=o(1).
$$
Since $\|u_k\|_{L^\infty(\O)}\le C$, by Lemma \ref{HHH}
we deduce that $$\io |\n (u-u_k)^+|^p \le o(1),$$  and, then, $$ \io |\n
(u_k-u)_-|^p \le o(1).$$ As a conclusion we get that
$$
\io |\n (u_k-u)_+|^p dx +\io |\n (u_k-u)_-|^p dx =o(1).
$$
Hence, $u_k\to u$ strongly in $W^{1,p}_0(\O)$ and the claim follows.\\

Finally, using the dominated convergence theorem, we deduce that
$$\dfrac{(u_+)_k^{\beta}|\nabla
u_k|^p}{1+\frac{1}{k}(u_+)_k^{\beta}|\nabla u_k|^p}\to u^\beta|\n
u|^p \mbox{  strongly in } L^1(\O).$$ It is clear that $u$ solves
the problem
\begin{equation}\label{apr01}
\left\{\begin{array}{rcll}
    -\Delta_p u +u^\beta|\nabla u|^p &=& f\;\; &\hbox{in } \Omega\\
    u &\geq & 0\;\; &\hbox{in } \Omega \\
    u &= & 0\;\;&\hbox{on } \partial\Omega.
\end{array}\right.
\end{equation}
By Proposition \ref{compagg}, the solution $u$ is unique. Moreover, letting $v=H(u)$, where
\begin{equation}\label{eq:fra1}
H(t)=\int_0^\s
e^{-\frac{t^{\beta+1}}{(p-1)(\beta+1)}}\, dt,
\end{equation}
we have that $v\in
W^{1,p}_0(\O),$
\begin{equation}\label{eq:28ap}
0\leq v\le \widehat{C}=\int_0^\infty
e^{-\frac{t^{\beta+1}}{(p-1)(\beta+1)}}\, dt,\end{equation}  and $v$
solves the problem
\begin{equation}\label{apr001}
\begin{cases}
-\De_p v=\dfrac{f}{(D'(v))^{p-1}}& \text{in } \Omega
\\
v=0 & \text{on } \partial \Omega,
\end{cases}
\end{equation}
where $D=H^{-1}$. By strong maximum principle, we deduce that  $v>0$ in $\Omega$ and from \eqref{eq:fra1}, it follows that  $u>0$ in $\Omega$ as well. Therefore, applying Proposition~\eqref{compagg}, we conclude  that  problem \eqref{apr01} admits   a unique positive bounded solution $u$. Moreover, if  $f_1\le f_2$, then the corresponding solutions satisfy  $u_1\le u_2$.

\

{\item \emph{Existence and uniqueness for $L^1$ datum.}} We assume now
that $f\in L^1(\O)$ is a nonnegative function. Define $u_n$ to be
the unique positive solution to the problem
\begin{equation}\label{apr01n}
\begin{cases}
-\De_p u_n+u_n^\beta|\nabla u_n|^p=T_n(f)& \text{in } \Omega
\\
u_n=0 & \text{on } \partial \Omega.
\end{cases}
\end{equation}
Using the first step we deduce that the sequence $\{u_n\}_n$ is
increasing in $n$. We claim that $\{u_n\}_n$ is bounded in
$W^{1,p}_0(\O)$. To prove the claim we use $T_k(u_n)$ as a test
function in \eqref{apr01n}, getting
$$
\io |\n T_k(u_n)|^p dx +\io T_k(u_n)u_n^\beta|\n u_n|^p dx \le
k\|f\|_{L^1(\O)}.
$$
Hence
$$
\io |\n T_k(u_n)|^p dx +k^{\beta+1} \int_{\{u_n\ge k\}}|\n u_n|^p
dx \le k\|f\|_{L^1(\O)}.
$$
Choosing $k=1$, we obtain that
\begin{equation}\label{eq:monit2}
\io |\n u_n|^p\le
\|f\|_{L^1(\O)}\end{equation} and\begin{equation}\label{pp}
\io T_k(u_n)u_n^\beta|\n u_n|^p dx \le k||f||_{L^1(\O)}.
\end{equation}
Therefore, we get the existence of $u\in W^{1,p}_0(\O)$ such that, up
to a subsequence, $u_n\rightharpoonup u$ weakly in $W^{1,p}_0(\O)$
and $u_n\rightarrow u$ a.e. in $\Omega$. \\
We claim that
$$\io u_n^\beta|\n u_n|^p dx \le C,$$ for all $n$. Indeed by \eqref{eq:monit2} and \eqref{pp} we obtain that
$$
\begin{array}{lll}
\dyle\io u_n^\beta|\n u_n|^p dx & = & \dyle\int_{u_n\le 1}
u_n^\beta|\n
u_n|^p dx +\int_{u_n\ge 1} u_n^\beta|\n u_n|^p dx\\
&\le & \dyle\io |\n u_n
|^p dx+\int_{u_n\ge 1} u_n^\beta|\n u_n|^p dx\\
&\le & \|f\|_{L^1(\O)},
\end{array}
$$
and, thus, the claim follows.

We now prove that $T_k(u_n)\to T_k(u)$ strongly in
$W^{1,p}_0(\O)$. To this end, we make use of  the change of variable
$v_n=H(u_n)$, as  defined above. It is clear that $v_n$ solves the
problem
\begin{equation*}
\begin{cases}
-\De_p v_n=\dfrac{f_n}{(D'(v_n))^{p-1}}& \text{in } \Omega
\\
v_n=0 & \text{on } \partial \Omega,
\end{cases}
\end{equation*}
where we set $f_n=T_n(f)$ and  $v_n\le  \widehat{C}$, see \eqref{eq:28ap}.   Moreover,  the sequence $\{v_n\}_n$ is
increasing in $n$. Observe that by using the definition \eqref{eq:28ap} once again  \eqref{eq:monit2}, the sequence  $\{v_n\}_n$ is bounded in
$W^{1,p}_0(\O)$. Since $-\D_p v_n\ge 0$, we can apply
Theorem~\ref{APM} to conclude  that $v_n \rightarrow v$ strongly in
$W^{1,p}_0(\O)$. Note also that for $k>0$, we have
$$
\{x\in \O: u_n(x){\leq} k\}=\{x\in \O: v_n(x){\leq}H(k)\}.
$$
By exploiting the strong convergence of the sequence $\{v_n\}_n$ in
$W^{1,p}_0(\O)$, we deduce that for every $a{\leq}H(\infty)=\widehat{C}$,
it holds that  $T_a(v_n)\to T_a(v)$ strongly in $W^{1,p}_0(\O)$. \\
Now, fixing  $k>0$ and defining  $a=H(k)$, clearly $a{\leq}\widehat{C}$. \\
Observe that
$$\n T_k(u_n)=\n u_n\chi_{\{u_n{\leq}k\}}=D'(v_n)\n
v_n\chi_{\{v_n{\leq}a\}}=D'(v_n)\n T_a (v_n).
$$
Since $D'(v_n)\le \widetilde{C}$ if $v{\leq}a{\leq}
\widehat{C}$, and applying the dominated convergence theorem, we
conclude  that
$$
D'(v_n)\n T_a(v_n)\to D'(v)\n T_a(v)\mbox{  strongly in
}(L^p(\O))^N.
$$
Noting that  $D'(v)\n T_a (v)=\n T_k(u)$, it follows  that $\n T_k(u_n)\to \n
T_k(u)$ strongly in $(L^p(\O))^N$. Thus, we conclude that $T_k(u_n)\to T_k(u)$
strongly in $W^{1,p}_0(\O)$.

Moreover, up to a subsequence, $\n u_n\to \n u$ a.e.  in $\O$.
To prove that $u_n\to u$ strongly in $W^{1,p}_0(\O)$ we apply
Vitali's Lemma. Let $E\subset \O$ be a measurable set. From
\eqref{eq:monit2} and \eqref{pp}, we obtain
\begin{equation}\label{eq:mang1}
\begin{array}{lll}
\dyle\int_E|\n u_n|^p dx & = & \dyle\int_{\{u_n\le k\}\cap E} |\n
u_n|^p dx +\int_{\{u_n\ge k\}\cap E} |\n u_n|^p dx\\
&\le & \dyle\int_E|\n T_k(u_n)
|^p dx+\frac{1}{k^\beta}\int_{u_n\ge k} u_n^\beta|\n u_n|^p dx\\
&\le & \dyle\int_E|\n T_k(u_n) |^p
dx+\frac{1}{k^\beta}||f||_{L^1(\O)}.
\end{array}
\end{equation}
Therefore,
$$\limsup_{n\to \infty}\int_E|\n u_n|^p dx \le \int_E|\n T_k(u)
|^p dx+\frac{C}{k^\beta}\le \int_E|\n u |^p dx+\frac{C}{k^\beta},
$$
where $k>0$ is fixed such that
${C}/{k^\beta}\le {\e}/{2}$. Since $u\in
W^{1,p}_0(\O)$, there exists  $\d>0$ such that for every measurable set $E\subset \Omega$ with  $|E|\le \d$, it holds that
\[\int_E|\n u |^p \,dx\leq \frac{\e}{2}.\]\\
Consequently, by  Vitali's Lemma we deduce that $|\nabla u_n|^p\rightarrow |\nabla u|^p$ in $L^1(\Omega)$.  Since  $W^{1,p}_0(\O)$ is  uniformly convex, this implies that
$u_n\to u$ strongly in
$W^{1,p}_0(\O)$.

To conclude this part, it remains  to prove that
$u_n^{\beta}|\nabla u_n|^p\to u^\beta|\n u|^p$ strongly in
$L^1(\O)$. Let  us define $G_k(s):=s-T_k(s)$, see \eqref{deftroncata}. Using $T_1(G_k(u_n))$  as a test function in
\eqref{apr01n}, one can show that
$$
\int_{u_n\ge k+1}u_n^\beta|\n u_n|^p dx \le \int_{u_n\ge k} f_n
dx.
$$
Now, applying  a similar argument to the one used in \eqref{eq:mang1},   we obtain  the equi-integrability of the sequence $\{u_n^\beta|\n u_n|^p\}_n$.
Hence, $u_n^{\beta}|\nabla u_n|^p\to u^\beta|\n u|^p$ strongly in
$L^1(\O)$.

Passing to the limit in \eqref{apr01n} we get the existence of a
solution $u$ to \eqref{one}.

\

The positivity  of $u$ follows from  the fact that the approximating sequence $\{u_n\}_n$ is  monotone increasing.

The solution $u$ is also unique: using the change of variable defined in \eqref{eq:28ap}, we have that $v=H(u)$ is a solution to
\begin{equation*}\begin{cases}
-\De_p v=\dfrac{f}{(D'(v))^{p-1}}& \text{in } \Omega
\\
v=0 & \text{on } \partial \Omega,
\end{cases}
\end{equation*}
where $D=H^{-1}$. This problem  admits a  unique solution, since  the term $1/(D'(\cdot))^{p-1}$ is decreasing and, thus, a weak comparison principle follows.
\end{enumerate}
\end{proof}
In the case where $f$ is replaced by a general Radon measure,
 the monotonicity of the sequence $\{u_n\}_n$ is no longer guaranteed. In
this situation, the limit problem may fail to have a  solution, which  typically occurs when the measure is singular. To clarify this point, we begin by defining what is meant by a  singular measure. We recall here below two classical definitions:
\begin{definition}\label{pcap}
Let $K$ be a compact subset of $\Omega$.
The $p$-capacity of $K$ with respect to $\Omega$ is defined by
\[
\operatorname{cap}_p(K) = \inf \left\{ \|\varphi\|_{W_0^{1,p}(\Omega)}^p : \varphi \in C_c^\infty(\Omega), \, \varphi \geq \chi_K \right\},
\]
where $\chi_K$ is the characteristic function of $K$
(and by convention we set  $\inf \varnothing = +\infty$).

The $p$-capacity of any open subset $U \subset \Omega$ is then defined by:
\[
\operatorname{cap}_p(U) = \sup \left\{ \operatorname{cap}_p(K) : K \text{ compact}, \; K \subset U \right\},
\]
and the $p$-capacity of any Borel set $B \subset \Omega$ by
\[
\operatorname{cap}_p(B) = \inf \left\{ \operatorname{cap}_p(U) : U \text{ open}, \; B \subset U \right\}.
\]
\end{definition}
\noindent Consider now a Radon measure $\mu$ on $\Omega$. We say
that $\mu$ is concentrated on a Borel subset $E\subseteq \Omega$, if
$\mu(B) = \mu(B\cap E)$ for every Borel set $B\subseteq\Omega$.
\begin{definition}\label{singular}
A positive Radon measure $\mu$ on $\Omega$ is said to be
{singular} if it is concentrated on a subset $E \subset \Omega$
such that $\operatorname{cap}_{p}(E)=0$.
We denote the set of singular positive measures by $\mathcal{M}_s$.\end{definition}
We now proceed with the

\

{\bf Proof of Theorem \ref{sing}.} Let $u_n$ be the unique positive solution to problem \eqref{apr01pp}
given in Theorem \ref{thone}. Since
$f_n\rightharpoonup \mu$, that is by definition,
\begin{equation}\label{eq:monit1}
\lim_{n\rightarrow+\infty}\int_{\Omega} f_n\varphi\, dx
=\int_\Omega \varphi \, d\mu
\quad\forall \varphi\in{C}(\overline{\Omega}),
\end{equation}
(see Remark 2.3 in \cite{DMOP}, for instance), we deduce that $\|f_n\|_{L^1(\O)}\le C$ (taking $\varphi=1$ in \eqref{eq:monit1}).  The above
estimates, (see equation \eqref{eq:monit2}), imply that $\{u_n\}_n$ is bounded in
$W_0^{1,p}(\Omega)$. Hence, there exists a  nonnegative
function $u\in W^{1,p}_0(\O)$ such that, up to subsequence,
$u_n\rightharpoonup u$ weakly in $W_0^{1,p}(\Omega)$.

Define $v_n=H(u_n)$, see \eqref{eq:fra1}. Then, $v_n\le \widehat{C}$ and $v_n$ solves
\eqref{apr001}. To show our result, we will prove that $v_n\to 0$
a.e. in $\O$.

Recall that $\mu\in \mathcal{M}_s$ is a singular Radon measure (see Definition \ref{singular}). Thus, for every $\epsilon>0$, there exists  a compact set $K_\epsilon\subset \Omega$ and   open sets $A_\epsilon,U_\epsilon\subset\subset \Omega$ such
that
\begin{equation}\label{eq:monit5}
K_\epsilon\subseteq E\subseteq A_\epsilon \quad \text{and}\quad\mu(A_\epsilon\setminus K_\epsilon)\leq \epsilon,\end{equation}
due to the regularity of the measure $\mu$ (Radon measures are regular measures), and moreover $K_\epsilon\subseteq U_\epsilon$.
By Definition \ref{pcap}, there exists $\varphi_\epsilon$
in ${C}_c^{\infty}(U_\epsilon)$ such that  $0\leq \varphi_\epsilon\leq1$ on
$U_\epsilon$, $\varphi_\epsilon\equiv 1$ on
$K_\epsilon$  and
$$ \int_{\Omega}|\nabla\varphi_\epsilon|^{p}dx<\epsilon.$$
Using Picone's inequality, we obtain
$$
\int_{\Omega} \dfrac{-\Delta_p v_n}{v^{p-1}_n}
\varphi_\epsilon^{p} dx\le
\int_{\Omega}|\nabla\varphi_\epsilon|^{p}dx\le \epsilon.
$$
Thus,
$$\int_{K_\epsilon}{\frac{f_n}{(v_n D'(v_n))^{p-1}}\,dx}\leq \epsilon.$$
Since $v_n\le \widehat{C}$, it follows that \begin{equation}\label{eq:monit4}\int_{K_\epsilon}{\frac{f_n}{(D'(v_n))^{p-1}}\,dx}=\int_{K_\epsilon}{\frac{f_n}{(v_n D'(v_n))^{p-1}}v_n^{p-1}\,dx}\leq
{C}^{p-1}\epsilon.\end{equation}

{We now claim that, for every $\phi\in
{C}^\infty(\overline{\Omega})$,
$$\lim\limits_{n\to
\infty}\int_{\Omega}{\frac{f_n}{(D'(v_n))^{p-1}}}\phi\,dx=0.$$ }
We decompose the integral into three parts:
\begin{eqnarray}\label{eq:monit3}\\\nonumber
&&\int_{\Omega}\frac{f_n}{(D'(v_n))^{p-1}}|\phi| \,dx =
\int_{K_\epsilon}\frac{f_n}{(D'(v_n))^{p-1}}|\phi| \,dx \\\nonumber
&&+\int_{A_\epsilon\setminus K_\epsilon}\frac{f_n}{(D'(v_n))^{p-1}}|\phi| \,dx +\int_{\Omega\setminus A_\epsilon}\frac{f_n}{(D'(v_n))^{p-1}}|\phi| \,dx\\\nonumber&&\le
\|\phi\|_\infty
\int_{K_\epsilon}\frac{f_n}{(D'(v_n))^{p-1}}\,dx+\int_{A_\epsilon\setminus
K_\epsilon}\frac{f_n}{(D'(v_n))^{p-1}}|\phi|\,dx+\int_{\Omega\setminus A_\epsilon}\frac{f_n}{(D'(v_n))^{p-1}}|\phi| \,dx
\\\nonumber
&&\leq  C\|\phi\|_\infty \epsilon+\int_{A_\epsilon\setminus
K_\epsilon}\frac{f_n}{(D'(v_n))^{p-1}}|\phi|\,dx+\int_{\Omega\setminus
A_\epsilon}\frac{f_n}{(D'(v_n))^{p-1}}|\phi|\,dx,
\end{eqnarray}
where we used \eqref{eq:monit4}.
Let us estimate the two last terms on the right-hand side of \eqref{eq:monit3}.
Since $(D'(v_n))^{1-p}$ is uniformly bounded, we have
\begin{equation}\label{eq:monit66}
\int_{A_\epsilon\setminus
K_\epsilon}\frac{f_n}{(D'(v_n))^{p-1}}|\phi|\,dx\leq C\int_{\Omega}f_n\chi_{A_\epsilon\setminus
K_\epsilon}|\phi|\,dx, \end{equation}
where $\chi_{A_\epsilon\setminus
K_\epsilon}$ denotes the  characteristic function of the set ${A_\epsilon\setminus
K_\epsilon}$.

Using \eqref{eq:monit1}, for every $\epsilon>0$,  there exists $n_0$ such that
\begin{eqnarray}\label{eq:monit6}
&&\int_{\Omega}f_n\chi_{A_\epsilon\setminus
K_\epsilon}|\phi|\,dx\leq\int_\Omega \chi_{A_\epsilon\setminus
K_\epsilon} |\varphi| \, d\mu +\epsilon\\\nonumber
&&=\int_{A_\epsilon\setminus
K_\epsilon} |\varphi| \, d\mu +\epsilon.
\end{eqnarray}
Therefore, by \eqref{eq:monit5}
\[\int_{A_\epsilon\setminus
K_\epsilon} |\varphi| \, d\mu \leq C\|\phi\|_\infty \epsilon.
\]
By \eqref{eq:monit66} and \eqref{eq:monit6} we obtain that
\begin{equation}\label{eq:monit7}
\int_{A_\epsilon\setminus
K_\epsilon}\frac{f_n}{(D'(v_n))^{p-1}}|\phi|\,dx\leq  C\|\phi\|_\infty \epsilon.
\end{equation}
Moreover, since $\mu(\Omega\setminus A_\epsilon)=0$, we may redefine $n_0$, if necessary, so that the sequence
 $\{f_n\}_n$ satisfies $\text{supp}\,(f_n)\subset
A_\epsilon$ for all $n\ge n_0$. Thus
\begin{equation}
\label{eq:monit8}
\int_{\Omega\setminus
A_\epsilon}{\frac{f_n}{(D'(v_n))^{p-1}}}|\phi|\,dx= 0\mbox{ if }
n\ge n_0.
\end{equation}
Hence, collecting \eqref{eq:monit3}, \eqref{eq:monit7} and  \eqref{eq:monit8}, we get
\begin{equation*}
\lim\limits_{n\to
\infty}\int_{\Omega}{\frac{f_n}{(D'(v_n))^{p-1}}}\,\phi\, \,dx=0,
\end{equation*}
and the claim follows.

Therefore, setting $\rho_n=\frac{f_n}{(D'(v_n))^{p-1}}$,  we obtain
$\rho_n\rightharpoonup 0$  in the narrow topology of measures (see  \eqref{eq:monit1}).

Since $\{u_n\}_n$ is bounded in $W^{1,p}_0(\O)$ and using the
fact that $H'\le 1$, we deduce that $\{v_n\}_n$ is bounded in
$W^{1,p}_0(\O)\cap L^\infty(\O)$ and that $v_n\rightharpoonup
v=H(u)$ weakly in $W^{1,p}_0(\O)$.

Recalling  that $-\D_p v_n=\rho_n$, then using Theorem 3.4 in
\cite{DMOP}, up to a subsequence, $\n v_n\to \n v$
a.e. in $\O$. Since for every  $\phi\in
{C}^\infty_c(\Omega)$,
$$
\lim\limits_{n\to \infty}\int_{\Omega}{-\Delta_p v_n\phi \,dx}=0,
$$
it follows that  $v\in W^{1,p}_0(\O)$ is a weak solution of
\begin{equation*}
-\Delta_p v=0\mbox{  in }\mathcal{D}'(\Omega).
\end{equation*}
Hence $v=0$, which implies $u=0$, see \eqref{eq:fra1}. This concludes the proof.
\begin{flushright}
    $\square$
\end{flushright}
\subsection{Analysis of   Problem \eqref{te}: Proof of Theorem~\ref{t1} and Theorem~\ref{singthh}}
In this section we prove Theorem \ref{t1}.
The solution that we are going to find is called {\em solution obtained as limits of approximations}, or simply SOLA, see \cite{DA} for further details.
By using the results in \cite{DMOP},  in this case SOLA is equivalent to an {\it entropy solution} or a {\it renormalized solution}.

The following well-known result will be useful.
\begin{Lemma}(Hardy-Sobolev Inequality)\label{le:Hardy}
Suppose  $1<p<N$ and $u\in W_0^{1,p}(\Omega)$. Then, we have
$$\int_{\mathbb{R}^N} \frac{|u|^p}{|x|^p}\leq C_{N,p}\int_{\mathbb{R}^N}|\nabla u|^p,$$
with ${\displaystyle C_{N,p}=\left(\frac{p}{N-p}\right)^p}$ optimal and not achieved constant.
\end{Lemma}
Now we are ready to prove the main result of the paper.

\

{\bf Proof of Theorem \ref{t1}: Case $(\ast)$}.
Let $u_0=0$. For $k\ge 1$, we define $u_k$ as the unique
positive solution (see Theorem \ref{thone} in the case of a bounded datum) of the approximating problem
\begin{equation}\label{eq:ta}
-\lap_p u_k+u_k^{\beta}|\nabla u_k|^p=\vartheta
T_k\left(\frac{u_{k-1}^q}{|x|^p}\right)+T_k(f)\;\;\hbox{in}\;\;\Omega,\;\;u_k\in
W_0^{1,p}(\Omega)\cap L^{\infty}(\Omega)
\end{equation}
From the result of  the previous section we deduce that $\{u_k\}_k$ is an
increasing sequence.\\
Taking $T_n(u_k)$ as a test function in the truncated problem
\eqref{eq:ta}, we obtain
\begin{eqnarray*}
\dyle \int_{\Omega}|\nabla
T_n(u_k)|^p+\int_{\Omega}u_k^{\beta}|\nabla u_k|^pT_n(u_k) & = &
\dyle \vartheta
\int_{\Omega}T_k\left(\frac{u_{k-1}^q}{|x|^p}\right)T_n(u_k)+\int_{\Omega}T_k(f)
T_n(u_k)\\
&\le & \dyle \vartheta
\int_{\Omega}T_k\left(\frac{u_{k}^q}{|x|^p}\right)T_n(u_k)+\int_{\Omega}T_k(f)
T_n(u_k).
\end{eqnarray*}

Notice that, for $n\ge 1$, we have
\begin{eqnarray*}
\int_{\Omega}u_k^{\beta}|\nabla
u_k|^p &= & \int_{u_k\le n}u_k^{\beta}|\nabla
u_k|^p+\int_{u_k>n}u_k^{\beta}|\nabla u_k|^p\\
&\le & n^\beta\int_{\O}|\nabla T_n(
u_k)|^p+\frac{1}{n}\int_{u_k>n}T_n(u_k)u_k^{\beta}|\nabla u_k|^p\\
&\le & n^\beta\bigg(\int_{\O}|\nabla T_n(
u_k)|^p+\int_{u_k>n}T_n(u_k)u_k^{\beta}|\nabla u_k|^p\bigg).
\end{eqnarray*}
Thus, for $n\ge 1$,
$$
\int_{\Omega}u_k^{\beta}|\nabla u_k|^p\le  \dyle n^{\beta+1}
\vartheta
\int_{\Omega}\frac{u_{k}^q}{|x|^p}+n^{\beta+1}||f||_{L^1(\O)}.
$$
Choosing $n=1$ and recalling that
$u_k^{\beta}|\nabla u_k|^p=
\left(\frac{p}{\beta+p}\right)^p\left |\n u_k^{\frac{\beta+p}{p}}\right|^p$, we obtain
\begin{equation}\label{madrid}
\left(\frac{p}{\beta+p}\right)^p\io |\n u_k^{\frac{\beta+p}{p}}|^p\le
\vartheta
\int_{\Omega}\frac{u_{k}^q}{|x|^p}+||f||_{L^1(\O)}.
\end{equation}
Using Hardy inequality (see Lemma \ref{le:Hardy}), we get
$$
\left(\frac{p}{\beta+p}\right)^p\io |\n u_k^{\frac{\beta+p}{p}}|^p\ge
\left(\frac{p}{\beta+p}\right)^p\left(\frac{N-p}{p}\right)^p\io
\frac{u_k^{p+\beta}}{|x|^p}.
$$
Hence
$$
\left(\frac{p}{\beta+p}\right)^p\left(\frac{N-p}{p}\right)^p\io
\frac{u_k^{p+\beta}}{|x|^p}\le \vartheta
\int_{\Omega}\frac{u_{k}^q}{|x|^p}+\|f\|_{L^1(\O)}.
$$
Since $q<p+\beta$,  Young's inequality yields
\begin{equation}\label{estim05}
\io \frac{u_k^{p+\beta}}{|x|^p}\le C_1+C_2||f||_{L^1(\O)},
\end{equation}
where $C_1,C_2$ are positive constants depending only on
$N,p,\beta$ and  independent of $u_k$. Therefore, the sequence
$\left\{\dfrac{u_k^{p+\beta}}{|x|^p}\right\}_k$ is uniformly
bounded in $L^1(\O)$ and so is the sequence
$\left\{T_k\left(\frac{u_{k-1}^q}{|x|^p}\right)\right\}_k$. Since
$\{u_k\}_k$ is increasing, if we set
\[h_k=\vartheta
T_k\left(\frac{u_{k-1}^q}{|x|^p}\right)+T_k(f),\] then $\{h_k\}_k$
is increasing and bounded in $L^1(\O)$. By the monotone
convergence theorem, there exists  $h\in L^1(\O)$ such
that $h_k\rightarrow h$ strongly in $L^1(\O)$. By Theorem \ref{thone}, there exists $u\in
W^{1,p}_0(\O)$ such that $u^\beta |\n u|^p\in L^1(\O)$, $u_k\to u$
strongly in $W^{1,p}_0(\O)$, $u_k^\beta |\n u_k|^p\to u^\beta |\n
u|^p$ strongly in $L^1(\O)$ and $u$ solves
$$
-\D_p u+u^\beta |\n u|^p=h.
$$
Recall that $\left\{\dfrac{u_{k-1}^{q}}{|x|^p}\right \}_k$ is
bounded in $L^1(\O)$. Applying once again the monotone convergence
theorem, we deduce that $\frac{u_{k-1}^{q}}{|x|^p}\to
\frac{u^{q}}{|x|^p}$ strongly in $L^1(\O)$. Therefore,
$h=\dfrac{u^{q}}{|x|^p}+f$ and the existence result follows.
Finally, we have $u>0$ in $\Omega$ since the sequence $\{u_k\}_k$
is monotone increasing.

\

{\bf Proof of Theorem \ref{t1}: Case $(\ast\ast)$.} Let us now set
$q=p+\beta$. Notice that, compared with the previous case, the main difficulty lies in the application of the Hardy inequality (see Lemma \ref{le:Hardy}). Recall
that $u_k$ is the solution to the truncated problem
\begin{equation*}
-\lap_p u_k+u_k^{\beta}|\nabla u_k|^p=\vartheta
T_k\left(\frac{u_{k-1}^{\beta
+p}}{|x|^p}\right)+T_k(f)\;\;\hbox{in}\;\;\Omega,\;\;u_k\in
W_0^{1,p}(\Omega)\cap L^{\infty}(\Omega),
\end{equation*}
where $\{u_k\}_k$ increasing in $k$. We claim that there exits
$\bar \vartheta=\overline\vartheta(p,\beta,N)$ such that,  for all
$0<\vartheta< \bar\vartheta$, the sequence $\left
\{\dfrac{u_k^{p+\beta}}{|x|^p}\right \}_n$ is bounded in
$L^1(\O)$.

To prove the claim, we take $T_n(u_k)$ as a test function in the
equation above. Recalling \eqref{madrid}, we obtain
\begin{equation*}
\left(\frac{p}{\beta+p}\right)^p\io |\n u_k^{\frac{\beta+p}{p}}|^p\le
\vartheta
\int_{\Omega}\left(\frac{u_{k}^{p+\beta}}{|x|^p}\right)+||f||_{L^1(\O)}.
\end{equation*}
Since, by Hardy inequality,
$$
\left (\frac{p}{\beta+p}\right)^p\io |\n u_k^{\frac{\beta+p}{p}}|^p\ge
\left (\frac{p}{\beta+p}\right)^p\left(\frac{N-p}{p}\right)^p\io
\frac{u_k^{p+\beta}}{|x|^p},
$$
we deduce that
$$
\left(\frac{p}{\beta+p}\right)^p\left(\frac{N-p}{p}\right)^p\io
\frac{u_k^{p+\beta}}{|x|^p}\le \vartheta
\int_{\Omega}\left(\frac{u_{k}^{p+\beta}}{|x|^p}\right)+||f||_{L^1(\O)}.
$$
Thus,  if
$\vartheta<\overline{\vartheta}:=\left (\frac{N-p}{{\beta+p}}\right)^p$,
we obtain
$$
(\overline{\vartheta}-\vartheta)\io \frac{u_k^{p+\beta}}{|x|^p}\le
||f||_{L^1(\O)}.
$$
Therefore, we conclude that the sequence
$\left\{\dfrac{u_k^{p+\beta}}{|x|^p}\right\}_k$ is bounded in $L^1(\O)$ and
the claim follows.

As a consequence of the compactness result in
Theorem \ref{thone} and using the same argument as in the previous case, we obtain the existence of the minimal solution $u$
to problem \eqref{te}.
\begin{flushright}
    $\square$
\end{flushright}
\begin{remark}
Define the set $$ \Upsilon=\{\vartheta>0\mbox{  such that
}\eqref{te} \mbox{  has a solution}\},
$$
then $\Upsilon$ is an interval.
\end{remark}
We now assume that
$f$ is a nonnegative measure which is singular with respect to the capacity $\operatorname{cap}_{p}(\Omega)$. We now give the

\

{\bf Proof of Theorem \ref{singthh}.} Notice that the existence of the minimal solution  $u_n$  to \eqref{sin0111} follows by  Theorem
\ref{t1}. In  its proof, it is shown that
$$
\io \frac{u_k^{p+\beta}}{|x|^p}\le C_1+C_2||f_n||_{L^1(\O)}\le C,
$$
see estimate \eqref{estim05}, where $\{u_k\}_k$ is the approximating sequence corresponding to the  minimal solution $u_n$ and where we have used the fact that $\|f_n\|_{L^1(\Omega)}\leq C$, by \eqref{eq:monit1}.

As  consequence, we deduce  that the
sequence $\{u_n\}_n$ is bounded in $W^{1,p}_0(\O)$, see Theorem~\ref{thone}. Hence, there exists  $\hat{u}\in W^{1,p}_0(\O)$ such that, up to a
subsequence, $u_n\rightharpoonup \hat{u}$ weakly in
$W^{1,p}_0(\O)$ and $u_n\to \hat{u}$ a.e. in $\O$. Setting
$v_n=H(u_n)$, where $H$ is defined in \eqref{eq:fra1}, we have  $v_n\le
\widehat{C}$ for all $n$ and $v_n$ solves the problem
\begin{equation*}
\begin{cases}
-\De_p v_n=
\vartheta\dfrac{D^q(v_n)}{(D'(v_n))^{p-1}|x|^p} + \dfrac{f_n}{(D'(v_n))^{p-1}}  & \text{in } \Omega
\\
v_n=0 & \text{on } \partial \Omega.
\end{cases}
\end{equation*}
Exactly as in the proof of Theorem \ref{sing}, we deduce that
$$\dfrac{f_n}{(D'(v_n))^{p-1}}\rightharpoonup 0$$ in the narrow
topology of measures. Since $v_n\le \widehat {C}$, we obtain
\[\vartheta\frac{D^q(v_n)}{(D'(v_n))^{p-1}|x|^p}\to
\vartheta\frac{D^q(\hat{v})}{(D'(\hat{v}))^{p-1}|x|^p}\quad \text{strongly in }L^1(\O).\] Thus, $\hat{v}\in  W^{1,p}_0(\O)\cap L^\infty(\O)$,
solves the problem
\begin{equation*}
\begin{cases}
-\De_p v=
\vartheta\dfrac{D^q(v)}{(D'(v))^{p-1}|x|^p}  & \text{in } \Omega
\\
v_n=0 & \text{on } \partial \Omega.
\end{cases}
\end{equation*}
Hence, $\hat{u}$ solves problem~\eqref{te} with $f=0$ and the
result follows.
\begin{flushright}
    $\square$
\end{flushright}
\section{Regularity and Behaviour at zero of the Solutions to \eqref{te}} \label{sec3}
We begin this section by proving that solutions to \eqref{te} are always unbounded, even in the case of a smooth source term $f$. This contrasts with problems without a Hardy potential, where, according to classical regularity theory, solutions corresponding to regular source terms are bounded. This unboundedness is in fact a consequence of the presence of the Hardy potential in the equation. Indeed, we prove something stronger: the solutions tend to infinity at zero (see Theorem \ref{l1}). This information is, for instance, fundamental in the study of qualitative properties such as symmetry and monotonicity of the solutions to problem~\eqref{te}.

\

\noindent{\bf Proof of Theorem \ref{l1}.} Consider the test function
\[\varphi=e^{-\frac{u^{\beta+1}}{\beta+1}}\psi,\] with $\psi \in W^{1,p}_0(\Omega)\cap L^{\infty}(\Omega)$ so that $\varphi \in W^{1,p}_0(\Omega)\cap L^{\infty}(\Omega)$. Using $\varphi$ as a test function in ~\eqref{te}, we obtain
\begin{equation}\label{eq:jgdasjhsgdfjhasghjiiiu}
\int_{\Omega}|e^{-\frac{u^{\beta+1}}{(\beta+1)(p-1)}} \nabla u|^{p-2}(e^{-\frac{u^{\beta+1}}{(\beta+1)(p-1)}}\nabla u, \nabla \psi)\geq\int_{\Omega}\frac{u^q}{|x|^p} e^{-\frac{u^{\beta+1}}{\beta+1}} \psi ,
\end{equation}
since $f$ is  nonnegative.  Recalling \eqref{eq:fra1}, from \eqref{eq:jgdasjhsgdfjhasghjiiiu}, we deduce   that $v=H(u)$ is a supersolution of
\begin{equation}\label{eq:hsajhasjhjsajhajhsaj1}\int_{\Omega}|\nabla v|^{p-2}(\nabla v, \nabla \psi)\geq\int_{\Omega}\frac{u^q}{|x|^p} e^{-\frac{u^{\beta+1}}{\beta+1}} \psi .\end{equation}

Now, let us consider  $u_R$ to be the radial solution to the problem
\begin{equation}\label{eq:gdajshgdhsjdgjhghdjsgrad}
\begin{cases}
-\Delta_p u + u^{\beta}|\nabla u|^p=\dfrac{C}{|x|^p} &\hbox{in}\; B_R,\\
u> 0 &\hbox{in}\;B_R,\\
u=0 &\hbox{on}\;\;\partial B_R.
\end{cases}
\end{equation}
The solution $u_R$ is radial, since by  Theorem \ref{thone} it is unique and $B_R$ is the standard ball of radius $R$.
Then, setting
\[\varphi=e^{-\frac{u_R^{\beta+1}}{\beta+1}}\psi\quad \text{and}\quad v_R=H(u_R),\]  we obtain
\begin{equation}\label{eq:lklksaasdfdf}
\int_{B_R}|\nabla v_R|^{p-2}(\nabla v_R, \nabla \psi)
=\int_{B_R}\frac{C}{|x|^p}e^{-\frac{u_R^{\beta+1}}{\beta+1}}\psi .
\end{equation}
By \eqref{eq:fra1}, we note that  $v$ (resp. $v_R$)  belongs to $W^{1,p}_0(B_R)\cap L^{\infty}(B_R)$ (to $W^{1,p}_0(B_R)\cap L^{\infty}(B_R)$). We may then  use in  \eqref{eq:hsajhasjhjsajhajhsaj1} and~\eqref{eq:lklksaasdfdf} the test function $\psi=(v_R-v)^+$,  where $R>0$ is chosen small enough so that $B_R \subset \subset \Omega$. In particular, since $v_R<v$  on~$\partial B_R$, we have $\psi \in W^{1,p}_0(\Omega)\cap L^{\infty}(\Omega)$. It follows that
\begin{equation}\label{eq:hsajhasjhjsajhajhsaj155}
\begin{split}
\int_{B_R}|\nabla v|^{p-2}(\nabla v, \nabla (v_R-v)^+) \,  \geq& \int_{B_R}\frac{u^q}{|x|^p}e^{-\frac{u^{\beta+1}}{\beta+1}}(v_R-v)^+ \,  \\
\geq&\int_{B_R}\frac{C_R}{|x|^p}e^{-\frac{u^{\beta+1}}{\beta+1}}(v_R-v)^+ \, ,
\end{split}
\end{equation}
with $C_R=\inf_{B_R}u^q(x)>0$ by the Harnack inequality (see \cite{PSB}). Moreover,
\begin{equation}\label{eq:lklksaasdfdf56}
\int_{B_R}|\nabla v_R|^{p-2}(\nabla v_R, \nabla(v_R-v)^+)
=\int_{B_R}\frac{C_R}{|x|^p}e^{-\frac{u_R^{\beta+1}}{\beta+1}}(v_R-v)^+ ,
\end{equation}
where, in \eqref{eq:gdajshgdhsjdgjhghdjsgrad}, we choose $C=C_R$.
Thus, subtracting
\eqref{eq:hsajhasjhjsajhajhsaj155} from  \eqref{eq:lklksaasdfdf56}, we obtain
$$\int_{B_R}(|\nabla v_R|^{p-2}\nabla v_R-|\nabla v|^{p-2} \nabla v, \nabla (v_R-v)^+)\leq 0.$$
Therefore, it follows that  $v\geq v_R$ in $B_R$. Using the definition of $v$ and $v_R$, we conclude that
\begin{equation}
\label{eq:shjahjdjhhdshdgsdghf111}
u  \geq u_R.
\end{equation}

Let us now study the qualitative behaviour of $u_R$.
First, observe  that the right hand side to \eqref{eq:lklksaasdfdf} is nonnegative, which implies  that $v_{R}$ is monotone decreasing  with respect to $\rho:=|x|$. Indeed, since $v_{R}$ is radial, for every $0<r\leq R$ we may use $(v_{R}(r)-v_{R}(x))^+\cdot \chi_{B_r}$ as a test function in \eqref{eq:lklksaasdfdf}, where $\chi_{B_r}$ denotes  the  characteristic function of the ball $B_r$. From this we deduce the desired monotonicity of $v_R$ with respect to  $\rho$.  Then, by  \eqref{eq:fra1}, we deduce that  $u_R(x)=u_R(|x|)$ is  also monotone decreasing  with respect to $\rho$ as well.

Consider   the test function  $\varphi=e^{\frac{-u_R^{\beta+1}}{\beta+1}}\psi$, with $\psi=\psi(|x|)$ belonging to $W^{1,p}_0(B_R)\cap L^{\infty}(B_R)$.  Then,  by \eqref{eq:gdajshgdhsjdgjhghdjsgrad},  we obtain
$$\int_0^Re^{-\frac{u_R^{\beta+1}}{\beta+1}}|u'_R|^{p-2}(u'_R, \psi')\rho^{N-1}d\rho=\int_{0}^RC_Re^{-\frac{u_R^{\beta+1}}{\beta+1}}\psi \rho^{N-1-p}d\rho.$$
By classical regularity results and Hopf's lemma, we have $u_R\in C^2(\overline B_R \setminus \{0\})$.  Hence,
$$\Big(e^{-\frac{u_R^{\beta+1}}{\beta+1}}|u_R'|^{p-2}(-u_R')\rho^{N-1}\Big)'= C_Re^{-\frac{u_R^{\beta+1}}{\beta+1}} \rho^{N-1-p} \quad \forall \, \rho \not= 0.$$
Since $u_R(\rho)$ is positive and monotone decreasing  with respect to $\rho$,  we distinguish between  the two following cases:

\

\begin{itemize}
\item[$(i)$] either $\displaystyle \lim_{\rho\rightarrow 0}u_R(\rho)=C>0$;

\

\item[$(ii)$] or $\displaystyle \lim_{\rho\rightarrow 0}u_R(\rho)=+\infty$.
\end{itemize}
If we assume  case~$(i)$, then
\[\frac{{\Big(e^{-\frac{u_R^{\beta+1}}{\beta+1}}|u_R'|^{p-2}(-u_R')\rho^{N-1}\Big)'}}{{(\rho^{N-p})'}}\rightarrow C \quad\text{as}\quad\rho \rightarrow 0,\]
for some positive constant $C$. It is standard to deduce that $-u_R'\geq{C}/{\rho}+o(1)$ for $\rho\rightarrow 0$,   getting a contradiction with the case $(i)$. Therefore,  case~$(ii)$ must hold, and together with \eqref{eq:shjahjdjhhdshdgsdghf111}, this concludes the proof.
\begin{flushright}
$\square$
\end{flushright}
Finally here below we prove our Calder\'on-Zygmund type result for the minimal solution to problem~\eqref{te}.

\

{\bf Proof of Theorem \ref{thm2bis}.} We prove the case $(*)$.

Let us
begin with assertion~$(i)$, namely we assume  that
$$q<p+\beta \qquad\text{and}\qquad m<\frac{N}{p}.$$
Fix $m$ as above, and let  $u_k$  be the unique solution to the approximating problem \eqref{eq:ta}. It is  known that the sequence $\{u_k\}_k\subset L^\infty(\Omega)$ is increasing in $k$, and that $u_k\to u$ strongly in $W^{1,p}_0(\O)$ where $u$ denotes  the minimal solution to problem \eqref{te}. Let $\a>0$ (to be chosen later). Using $u^{\a}_k$ as a test function in \eqref{eq:ta}, we obtain
\begin{equation}\label{eq:bell1}
    \begin{split}
        \int_{\Omega} u^{\alpha+\beta}_k  |\nabla u_k|^p \ \leq \vartheta \int_{\Omega} \frac{u^{q+\alpha}_k}{|x|^p}  + \int_\Omega f_k u^{\alpha}_k.
    \end{split}
\end{equation}
The function $u_k^\alpha$ can be used as test function in \eqref{eq:ta} by means of a standard approximation argument. Specifically,  we use the test function $(\delta +u_k)^{\alpha}-\delta^\alpha$, with $\delta>0$ in \eqref{eq:ta}, and then  let  $\delta\rightarrow 0$, which yields \eqref{eq:bell1}.
 Since $q<p+\beta$, the  weighted  Young  inequality together with Lemma~\ref{le:Hardy} implies  that, for every $\e>0$,
$$
\int_{\Omega} \frac{u^{q+\alpha}_k}{|x|^p} \le \e \vartheta\int_{\Omega} u^{\alpha+\beta}_k  |\nabla u_k|^p+ C(p,q,\beta,\e, \Omega).
$$
Using H\"older inequality, we also have
$$
\int_\Omega f_k u^{\alpha}_k\le
\left(\int_{\Omega}f^m \right)^{\frac 1m} \left(\int_{\Omega}u^{\alpha m'}_k \right)^{\frac{1}{m'}}\le C
\left(\int_{\Omega}u^{\alpha m'}_k \right)^{\frac{1}{m'}}.
$$
Hence,
\begin{equation}\label{roma}
        \int_{\Omega} u^{\alpha+\beta}_k  |\nabla u_k|^p \ \leq  C
        \left(\int_{\Omega}u^{\alpha m'}_k \right)^{\frac{1}{m'}}+ C.
\end{equation}
Choose
$$\a=\frac{p^*(m-1)(\beta+p)}{pm-(m-1)p^*},$$ so that
$\a m'=(\beta+p+\a)\frac{p^*}{p}$.
Note that $\alpha >0$, since $m<\frac Np$. By the  Sobolev inequality, we then obtain
$$
C_1\left(\int_{\Omega}u^{\alpha m'}_k \right)^{\frac{p}{p^*}}\le C_2\left(\int_{\Omega}u^{\alpha m'}_k \right)^{\frac{1}{m'}}+C_3,
$$
where $C_i>0$ are independent of $u_k$. Since $\frac{1}{m'}<\frac{p}{p^*}$, (because  $m<\frac Np$), we conclude that
$$
\int_{\Omega}u^{\alpha m'}_k\le C,\quad \mbox{for all }k>0.
$$
Therefore,  by Fatou's lemma and setting $\gamma =\alpha m'$, we obtain the conclusion.

Returning to \eqref{roma} and applying Fatou's  lemma once more, we deduce that
$$\int_{\Omega} u^{\alpha+\beta}|\nabla u|^p\leq  C,$$
which completes the proof of case  $(i)$, noting that
\[\alpha+\beta= {\frac{p(m\beta+p^*(m-1))}{pm-p^*(m-1)}}.\]

We  now   prove  case $(ii)$. We assume that
\[q<p+\beta \qquad\text{and}\qquad m\geq\frac{N}{p}.\]
It is clear that $f\in L^{m_1}(\O)$ for all
$m_1<\frac{N}{p}$. Setting $\e=pm_1-p^*(m_1-1)$, by the previous  case we obtain, for all $\e>0$,
$$
\io u^{\frac{p(\beta m_1+p^*(m_1-1))}{\e}}|\n u|^p\le C.
$$
Observing that $$\frac{p(\beta m_1+p^*(m_1-1))}{\e}\rightarrow +\infty \quad \text{as}\quad\varepsilon \rightarrow 0,$$ (namely, when $m_1 \rightarrow \frac Np$), we obtain the desired conclusion for case $(ii)$.

\

Finally,  case $(**)$ follows by repeating the same argument used for case $(*)$ together with  Lemma \ref{le:Hardy}.
\begin{flushright}
    $\square$
\end{flushright}
{\bf Comments and remarks.} As a direct and relevant
consequence of Theorems \ref{l1} and ~\ref{thm2bis}, we deduce  that, if $f\in L^m(\Omega)$ with $m>\frac{N}{p}$, the corresponding solution    $u$ belongs to $L^\gamma(\Omega)$ for every $\gamma< +\infty$. Nevertheless, as already observed at the beginning of this section,  Theorem \ref{l1} ensures that   $u\notin L^\infty(\O)$.

In certain specific cases, we can describe more precisely the behaviour of the solution near the singularity at the origin.  We only wish to provide an idea and a sketch of the proof, since the argument is rather  technical. Given $R>0$, let us set \[w(x)=\ln\left(\frac{R}{|x|}\right).\] Clearly,  $w\in W^{1,p}_0(B_R)$, and   for $|x|<\eta<R$, we have
\begin{equation*}
-\lap_p w+w^{\beta}|\nabla w|^p= \frac{N-p}{|x|^p}+\frac{1}{|x|^p}\left(\ln\left(\frac{R}{|x|}\right)\right)^\beta.
\end{equation*}
\begin{enumerate}
\item Suppose that $\beta\geq q$ and
$$f(x)=\frac{1}{|x|^\theta},$$
with $\theta<p$. In this case $f\in L^m(\Omega)$ for some  $m>\frac{N}{p}$.  There exists a   constant $C>0$ such that  $w_1=Cw$ satisfies
$$
-\lap_p w_1+w_1^{\beta}|\nabla w_1|^p\ge \frac{w_1^q}{|x|^p}+f,$$
for $x\approx 0$. Hence, if $u$ denotes  the minimal solution of problem
\eqref{te} with $f(x)=\frac{1}{|x|^\theta}$,  then assuming $R$ large enough so that $\Omega\subset B_R$,  the minimality property of $u$ (see Remark \ref{rem:min}, $(i)$) combined with  the  comparison principle in Proposition \ref{compagg}
 yields  $$u(x)\le C\ln\left(\frac{R}{|x|}\right),$$ for $x$ near
the origin.
\item Let $\beta<p-1<q$, and define \[f(x)=\frac{1}{|x|^p}.\] Let $u$ be the minimal
solution of problem \eqref{te}.
For some $\eta>0$, we have
\begin{equation}\label{eq:comp1}
-\lap_p w_1+w_1^{\beta}|\nabla w|^p\le  \frac{w_1^\beta}{|x|^p}+f \mbox{  in   }B_\eta(0),
\end{equation}
and
\begin{equation}\label{eq:comp2}
-\lap_p u+u^{\beta}|\nabla u|^p\ge  \frac{u^\beta}{|x|^p}+f \mbox{  in   }B_\eta(0),
\end{equation}
where we have used the fact that, thanks to Theorem \ref{l1}, one has $u(x)\to \infty$ as $|x|\to 0$.
Using the change of variable \eqref{eq:fra1}, inequalities \eqref{eq:comp1} and \eqref{eq:comp2} become
\[ -\lap_p v_1\leq\frac{1}{|x|^p} \frac{D(v_1)^\beta}{(D'(v_1))^{p-1}}\mbox{  in   }B_\eta(0),
 \]
 and
 \[
 -\lap_p v\geq\frac{1}{|x|^p} \frac{D(v)^\beta}{(D'(v))^{p-1}}, \mbox{  in   }B_\eta(0),
\]
where $v_1=H(w_1)$ and $v=H(u)$ and $D=H^{-1}$.

Now take $R=\eta$. By  the explicit definition of $w$, it follows that $$K:=\text{supp}(w_1-u)^+\subset B_\eta.$$ Therefore, $w,u\geq \delta>0$ in $K$ and by \eqref{eq:fra1} it holds $v_1,v\geq \delta'>0$. Moreover, we have  $v_1, v\in W^{1,p}_0(\Omega)\cap L^{\infty}(\Omega)$ such that
\[\frac{v_1}{v},\,\, \frac{v}{v_1}\in L^{\infty}(\Omega).\]
Thus, Proposition \ref{rome2} applies, yielding $v\geq v_1$ in $B_\eta(0)$. In terms of the original function $u,w_1$, this implies
$$u(x)\ge C\ln\left(\frac{1}{|x|}\right),$$
for $x$ close
to the origin.
\end{enumerate}
\begin{remark}
It appears natural to conjecture that, in the case $q=p+\beta$,
the condition~$(\ast\ast)$ stated in Theorem~\ref{t1} is indeed
optimal. More precisely, we expect than when the parameter $\vartheta$ is
sufficiently large,  problem \eqref{te} {does
not admit any positive solution}.
\end{remark}

\begin{center}
{\sc Acknowledgements}
\end{center} 
B.  Abdellaoui is partially supported by the Project PID2023-150166NB-I00,
{\em EDPs no-lineales: flujos, desigualdades y geometr\'ia}.

L. Montoro is partially supported by the PRIN project P2022YFAJH 003 (Italy), {\em Linear and nonlinear PDEs; new directions
and applications}, and  by the Gruppo Nazionale per l'Analisi Matematica, la Probabilit\`a  e le loro Applicazioni (GNAMPA) of the Istituto Nazionale di Alta Matematica (INdAM).

\vspace{1cm}

\

\begin{center}
{\sc Data availability statement}

\

All data generated or analyzed during this study are included in this published article.
\end{center}

\

\begin{center}
	{\sc Conflict of interest statement}
	
	\
	
	The authors declare that they have no competing interest.
\end{center}


\end{document}